\documentclass[reqno,centertags]{amsart}
\usepackage{amsthm,amsmath,amsfonts,bm,amssymb,epsfig}
\usepackage{mathptmx}
\usepackage{mathrsfs}

\DeclareSymbolFont{SY}{U}{psy}{m}{n}
\DeclareMathSymbol{\emptyset}{\mathord}{SY}{'306}

\newtheorem{introtheorem}{Theorem}{\bf}{\it}

\newtheorem{theorem}{Theorem}[section]

\newtheorem{lemma}[theorem]{Lemma}

\newtheorem{hypothesis}[theorem]{Hypothesis}
\theoremstyle{definition}

\theoremstyle{remark}
\newtheorem{introremark}[introtheorem]{Remark}{\it}{\rm}
\newtheorem{remark}[theorem]{Remark}
\newtheorem{example}[theorem]{Example}
\numberwithin{equation}{section}

\newcommand{\bbC}{\mathbb{C}}

\newcommand{\bbR}{\mathbb{R}}

\newcommand{\sE}{\mathsf{E}}

\newcommand{\fH}{\mathfrak{H}}
\newcommand{\fK}{\mathfrak{K}}

\newcommand{\lal}{\langle}
\newcommand{\ral}{\rangle}
\newcommand{\got}{\mathfrak}

\newcommand{\ran}{\mathop{\mathrm{Ran}}}
\newcommand{\Ran}{\mathop{\mathrm{Ran}}}

\newcommand{\spec}{\mathop{\text{\rm spec}}}
\newcommand{\dist}{\mathop{\mathrm{dist}}}
\newcommand{\dom}{\mathop{\mathrm{Dom}}}
\newcommand{\Dom}{\mathop{\mathrm{Dom}}}

\begin{document}

\title[The \textit{a priori} tan\,$\Theta$ Theorem for Eigenvectors]
{The \textit{a priori} tan\,$\Theta$ Theorem for Eigenvectors}

\author[S. Albeverio, A. K. Motovilov, and A. V. Selin]
{S.\,Albeverio, A.\,K.\,Motovilov, and A.V. Selin}

\address{Sergio Albeverio,
Institut f\"ur Angewandte Mathematik, Universit\"at Bonn,
Wegelerstra{\ss}e 6, D-53115 Bonn, Germany; SFB 611, Bonn; BiBoS,
Bielefeld-Bonn; CERFIM, Locarno; Accademia di Architettura, USI, Mendrisio}
\email{albeverio@uni-bonn.de}

\address{Alexander K. Motovilov, Bogoliubov Laboratory of
Theoretical Physics, JINR, Joliot-Curie 6, 141980 Dubna, Moscow
Region, Russia} \email{motovilv@theor.jinr.ru}

\address{Alexei V. Selin, Open Technologies Inc., Obrucheva 30 str. 1,
117997 Moscow, Russia} \email{selin@ot.ru}

\subjclass{Primary 47A55; Secondary 47B25}

\keywords{Perturbation problem, spectral subspaces, perturbation of eigenvectors,
$\tan\theta$ theorem}

\date{December 22, 2005}
\begin{abstract}
Let $A$ be a self-adjoint operator on a Hilbert space $\fH$. Assume
that the spectrum of $A$ consists of two disjoint components
$\sigma_0$ and $\sigma_1$ such that the convex hull of the set $\sigma_0$
does not intersect the set $\sigma_1$. Let $V$ be a bounded
self-adjoint operator on $\fH$ off-diagonal with respect to the
orthogonal decomposition $\fH=\fH_0\oplus\fH_1$ where $\fH_0$ and
$\fH_1$ are the spectral subspaces of $A$ associated with the
spectral sets $\sigma_0$ and $\sigma_1$, respectively. It is known
that if $\|V\|<\sqrt{2}d$ where $d=\dist(\sigma_0,\sigma_1)>0$ then
the perturbation $V$ does not close the gaps between $\sigma_0$ and
$\sigma_1$. Assuming that $f$ is an eigenvector of the perturbed
operator $A+V$ associated with its eigenvalue in the interval
$(\min(\sigma_0)-d,\max(\sigma_0)+d)$ we prove that under the
condition $\|V\|<\sqrt{2}d$ the (acute) angle $\theta$ between $f$
and the orthogonal projection of $f$ onto $\fH_0$ satisfies the
bound $\tan\theta\leq\frac{\|V\|}{d}$ and this bound is sharp.

\end{abstract}

\maketitle

\section{Introduction}

Given a self-adjoint operator $A$ on a Hilbert space
$\fH$ assume that $\sigma_0$ is an isolated part
of its spectrum, that is,
\begin{equation}
\label{ddist}
d=\dist(\sigma_{0}, \sigma_{1})  > 0,
\end{equation}
where $\sigma_1=\spec(A)\setminus\sigma_0$ is the rest of the
spectrum of $A$. In this case we say that there are open gaps
between the sets $\sigma_0$ and $\sigma_1$. It is well known (see,
e.g., \cite[\S 135]{SzNagy}) that a sufficiently small self-adjoint
perturbation $V$ of $A$ does not close these gaps which allows one
to think of the corresponding disjoint spectral components
$\sigma'_{0}$ and $\sigma'_{1}$ of the perturbed operator $L=A + V$
as a result of the perturbation of the spectral sets $\sigma_{0}$
and $\sigma_{1}$, respectively.

Assuming \eqref{ddist}, in this note we are concerned with the
perturbations $V$ that are off-diagonal with respect to the
partition $\spec(A)=\sigma_0\cup\sigma_1$, i.e. with perturbations
that anticommute with the difference
$\sE_A(\sigma_{0})-\sE_A(\sigma_{1})$ of the spectral projections
$\sE_A(\sigma_{0})$ and $\sE_A(\sigma_{1})$ associated with the
spectral sets $\sigma_0$ and $\sigma_1$, respectively. In general,
it is known (see \cite[Theorem 1]{KMM4}) that such perturbations do
not close the gaps between the sets $\sigma_0$ and $\sigma_1$ (which
means that the inequality $\dist(\sigma'_0,\sigma'_1)>0$ holds)
whenever
\begin{equation}
\label{sqrt32} \| V \| < \frac{\sqrt{3}}{2} d.
\end{equation}
Moreover, if no assumptions are made about the location of
$\sigma_0$ and $\sigma_1$ except the assumption \eqref{ddist}
then condition \eqref{sqrt32} is sharp (see \cite[Example
1.5]{KMM4}).

However there are two important particular mutual dispositions
of the spectral sets $\sigma_0$ and $\sigma_1$ that ensure the
disjointness of the perturbed spectral sets $\sigma'_0$ and
$\sigma'_1$ under conditions on $\|V\|$ much weaker than the
general one of \eqref{sqrt32}. The first of these two
dispositions is the one where the sets $\sigma_0$ and $\sigma_1$
are subordinated, say
\begin{equation}
\label{tan2t} \sup \sigma_{0} <  \inf \sigma_{1}.
\end{equation}
The second disposition corresponds to the case where one of the
sets $\sigma_0$ and $\sigma_1$ is lying in a finite gap of the other
set, say $\sigma_0$ lyes in a finite gap of $\sigma_1$, which means that
\begin{equation}
\label{tant}
\mathop{\mathrm{conv}}( \sigma_{0})\cap\sigma_{1}
= \emptyset,
\end{equation}
where $\mathop{\mathrm{conv}}(\sigma)$ denotes the convex hull
of a set $\sigma\subset\bbR$. (We recall that by a finite gap of a
closed Borel set $\Sigma$ on $\bbR$ one understands an open finite
interval belonging to the complement $\bbR\setminus\Sigma$ of $\Sigma$
such that both of its end points belong to $\Sigma$.)

It is known that if \eqref{tan2t} holds then for any bounded
off-diagonal perturbation $V$ the interval
$(\sup\sigma_0,\inf\sigma_1)$ belongs to the resolvent set of
the perturbed operator $L=A+V$, and hence $\sigma'_{0} \subset
(-\infty, \sup \sigma_{0}]$ and $\sigma'_{1} \subset [\inf
\sigma_{1}, +\infty)$ (see \cite{ALT}, \cite{DK}; cf.
\cite{KMM3}). In the case of the disposition \eqref{tant}, it
has been proven in \cite{KMM4} (see also \cite{KMM3}) that the
gaps between $\sigma_0$ and $\sigma_1$ remain open if the
off-diagonal perturbation  $V$ satisfies the (sharp) condition
\begin{equation*}
\| V \| < \sqrt{2} d.
\end{equation*}
Under this condition the spectrum of $L=A+V$ consists of two
disjoint components $\sigma'_0$ and $\sigma'_1$ such that
$$
\sigma_{0}'\subset(\inf \, \sigma_{0} -
d, \sup \, \sigma_{0} + d )
\text{\, and \,}
\sigma'_{1}\subset\bbR\setminus\Delta,
$$
where $\Delta$ denotes the gap of $\sigma_1$ that contains
$\sigma_0$. Notice that the norm bound $\|V\|<\sqrt{2}d$ is also
sharp in the sense that, if it is violated, the spectrum of $L$ in
the gap $\Delta$ may be empty  at all (see \cite[Example
1.6]{KMM4}).

Now assume that the perturbed spectral set $\sigma'_0$ contains
an eigenvalue of the operator $L=A+V$ and let $f$, $f\neq 0$, be an
eigenvector of $L$ corresponding to this eigenvalue. Denote by
$\theta$ the (acute) angle between the vector $f$ and its
projection $f_0=\sE_A(\sigma_0)f$ onto the spectral subspace
$\fH_0=\Ran\sE_A(\sigma_0)$ of $A$ associated with the
unperturbed spectral set~$\sigma_0$.

It is known that under the subordination condition \eqref{tan2t}
for any bounded off-diagonal perturbation $V$ the angle $\theta$ can not
exceed ${\pi}/{4}$. Moreover, the following sharp estimate
holds
\begin{equation}
\label{tan2t-a}
\theta\leq\frac{1}{2}\arctan\left(\frac{2\|V\|}{d}\right)\quad
\left(<\frac{\pi}{4}\right).
\end{equation}
This bound is a simple corollary to the celebrated Davis--Kahan
$\tan 2\Theta$ Theorem \cite{DK} (also see \cite[Theorem 6.1]{Davis:123}
and \cite[Theorem 2.4]{KMM5}).

In the case of the spectral disposition
\eqref{tant} an \textit{a posteriori} bound on the angle $\theta$
under condition $\|V\|<\sqrt{2}d$ follows from \cite[Theorem 2.4]{KMM4}.
This bound reads
\begin{equation}
\label{tandel}
\theta\leq\arctan\left(\dfrac{\|V\|}{\delta}\right),
\end{equation}
where $\delta$ denotes the distance between the perturbed spectral
set $\sigma'_0$ and unperturbed spectral set $\sigma_1$. Since
$\delta$ may be arbitrarily small (see Example \ref{remdel} below),
the bound \eqref{tandel} gives in general no \textit{a priori}
uniform estimate for $\theta$ except that $\theta<{\pi}/{2}$.

The present note is aimed just at giving an \textit{a priori} sharp bound
on the angle $\theta$ in the case of the disposition
\eqref{tant}. In particular, we will prove that under
condition $\|V\|<\sqrt{2}d$ this angle is strictly separated from
$\pi/2$. Our main result is as follows

\begin{introtheorem}
\label{Th1}
Given a self-adjoint operator $A$ on the Hilbert
space $\got H$ assume that
\begin{equation*}
\spec(A) = \sigma_{0}\cup\sigma_{1},  \quad
\dist(\sigma_{0},\sigma_{1}) = d > 0, \quad  \text{and} \quad
\mathop{\mathrm{conv}}(\sigma_0)\cap\sigma_1=
\emptyset.
\end{equation*}
Let $V$ be a bounded self-adjoint operator on $\fH$
off-diagonal with respect to the decomposition $\got H
= \ran \sE_{A}(\sigma_{0})\oplus\ran
\sE_A(\sigma_{1})$. Assume in addition that
\begin{equation}
\label{Vd2}
\| V \| < \sqrt{2}d
\end{equation}
and that the operator $L=A+V$ possesses an eigenvector $f$
associated with an eigenvalue
$$
z\in(\inf\,\sigma_{0}-d,\sup\,\sigma_{0}+d).
$$
Then the (acute) angle $\theta$ between the vector $f$ and its
projection \,$\sE_A(\sigma_0)f$\, onto the subspace
$\Ran\sE_A(\sigma_0)$ satisfies the bound
\begin{equation}
\label{thet}
\theta\leq\arctan\left(\frac{\|V\|}{d}\right).
\end{equation}
\end{introtheorem}
\begin{introremark}
The bound \eqref{thet} implies that under condition \eqref{Vd2} the angle
$\theta$ can never exceed the value of $\arctan\sqrt{2}$, i.e.
$$
\theta<\arctan\sqrt{2}\approx 0.304\,\pi.
$$
\end{introremark}
We also remark that for $\|V\|<d$ the bound \eqref{thet} follows
from \cite[Theorem 2.4]{MS01}.

Throughout the paper by $\Xi(D,d,b)$ we will denote a function
of three real variables $D$, $d$, and $b$ defined on the set
$$
\Omega=\bigl\{(D,d,b)\,\,|\quad D>0,\quad 0<d\leq D/2,\quad 0\leq
b<\sqrt{dD}\bigr\}
$$
by the following expressions
\begin{equation}
\label{Xi} \Xi(D,d,b)=\left\{
\begin{array}{l}
\tan^2\biggl(\dfrac{1}{2}\arctan\dfrac{2b}{d}\biggr) \quad\text{if
\,}
b^2\leq d\sqrt{D}\dfrac{\sqrt{D}-\sqrt{2d}}{2},\\[4mm]
1+\dfrac{2b^2}{D^2}- \dfrac{2}{D^2}\sqrt{
(dD-b^2)\bigl((D-d)D-b^2\bigr)} \\
\qquad\qquad\qquad\qquad\qquad \text{if \,\,}
{d\sqrt{D}}\dfrac{\sqrt{D}-\sqrt{2d}}{2}<b^2<dD
\end{array}
    \right.
\end{equation}
Here and further on by $\tan^2\theta$, $\theta\in\bbR$, we
understand the square of the tangent of $\theta$, that is,
$\tan^2\theta=(\tan\theta)^2$.

Theorem \ref{Th1} appears to be a corollary to a more
general statement (Theorem \ref{ThMain}) that is proven
under a weaker than \eqref{Vd2} but more specific
condition $\|V\|<\sqrt{d|\Delta|}$ where $\Delta$ again
denotes the (finite) gap of the set $\sigma_1$ that
contains $\sigma_0$ and $|\Delta|$ stands for the
length of the interval $\Delta$. If this condition
holds then the off-diagonal perturbation $V$ does not close
the gaps between $\sigma_0$ and $\sigma_1$ (see
\cite[Theorem 1 (i)]{KMM3}). The claim of Theorem
\ref{ThMain} is that under the condition
$\|V\|<\sqrt{d|\Delta|}$ the following inequality holds
\begin{equation}
\label{bo2in}
\tan\theta\leq\bigl(\Xi(|\Delta|,d,\|V\|)\bigr)^{1/2}.
\end{equation}
In particular, from formula \eqref{Xi} defining the function $\Xi$
one can see that if $|\Delta|>2d$ then for $V$ small enough,
namely for $V$ such that
$$
\|V\|^2\leq d\sqrt{|\Delta|} \dfrac{\sqrt{|\Delta|}-\sqrt{2d}}{2},
$$
the bound on $\theta$ is the same as the bound \eqref{tan2t-a}
prescribed by the $\tan2\Theta$ Theorem.

The paper is organized as follows. In Section \ref{Beginning} we
consider a three-dimensional version of the problem and prove the
bound \eqref{bo2in} in the case of $3\times3$ matrices. The general
infinite-dimensional case is studied in Section \ref{Sgener}. In the
proof of the central result of this section, the one of
Theorem~\ref{ThMain}, we essentially rely on Lemma \ref{lrot1} of
Section \ref{Beginning}.

Throughout the paper we use the standard notation $M^\intercal$
for the transpose of a matrix $M$.

\section{A three-dimensional case}
\label{Beginning}

We start our consideration with the case where $\fH=\bbC^3$ and the
operators $A$ and $V$ are $3\times3$ matrices. Assume that
\begin{equation*}
A=\left(\begin{array}{rll}
\lambda & 0 & 0 \\
0 & \gamma_- & 0 \\
0 &  0 & \gamma_+
\end{array}\right)
\text{\, and \,} V=\left(\begin{array}{lll}
0 & b_- & b_+ \\
b_- &  0 & 0 \\
b_+ &  0 & 0
\end{array}\right),
\end{equation*}
where
\begin{equation*}
\lambda,\gamma_\pm,b_\pm\in\bbR,  \text{\, and \,}
\gamma_+>\gamma_-.
\end{equation*}
The matrices $A$ and $V$ are symmetric. Moreover, under the
assumption that $\lambda\neq\gamma_\pm$ the matrix $V$ is
off-diagonal with respect to the partition
$\spec(A)=\sigma_0\cup\sigma_1$ of the spectrum of $A$ into the
disjoint sets
\begin{equation*}
\sigma_0=\{\lambda\} \text{\, and \,}
\sigma_1=\{\gamma_-,\gamma_+\}.
\end{equation*}
It is convenient for us to write the matrix $L=A+V$ in the following
$2\times2$ block form
\begin{equation}
\label{Lbl} L=\left(\begin{array}{cc}
\lambda & B \\
B^*    &  A_1
\end{array}\right),
\end{equation}
where $B$ and $A_1$ are $1\times 2$ and $2\times2$ matrices given by
\begin{equation}
\label{BA1}
 B=(\,b_-\quad b_+\,), \quad
A_1=\left(\begin{array}{rl} \gamma_- & 0 \\ 0 &
\gamma_+\end{array}\right),
\end{equation}
respectively. Clearly, $\|V\|=\|B\|=\sqrt{|b_-|^2+|b_+|^2}$.

Throughout this section by $\Delta$ we will denote the spectral gap
of the operator $A_1$ between its eigenvalues $\gamma_-$ and
$\gamma_+$, i.e.
$$
\Delta=(\gamma_-,\gamma_+).
$$

\begin{lemma}
\label{root3} Given a matrix $L$ of the form \eqref{Lbl},
\eqref{BA1}, assume that $\lambda\in\Delta$ and
\begin{equation}
\label{Bdd}
\|B\|<\sqrt{d|\Delta|},
\end{equation}
where $|\Delta|=\gamma_+ -\gamma_-$ stands for the length of the
interval $\Delta$ and\textbf{\emph{ }}
$d=\dist(\sigma_0,\sigma_1)=\min\{\gamma_+-\lambda,\lambda-\gamma_-\}$.
Then $L$ has a unique eigenvalue $z$ in the interval $\Delta$ and
this eigenvalue is simple. Moreover,
\begin{equation*}
\gamma_-<z_{\mathrm{min}}\leq z \leq  z_{\mathrm{max}}<\gamma_+,
\end{equation*}
where
\begin{align}
\label{dBm}
z_{\mathrm{min}}&=\lambda-\|B\|\tan\left(\frac{1}{2}\arctan\frac{2\|B\|}
{\gamma_+ -\lambda}\right),\\
\label{dBp}
z_{\mathrm{max}}&=\lambda+\|B\|\tan\left(\frac{1}{2}\arctan\frac{2\|B\|}
{\lambda-\gamma_-}\right).
\end{align}
\end{lemma}

\begin{proof} Lemma \ref{root3} is an elementary
corollary to \cite[Theorem 3.2]{KMM3}.
\end{proof}

\begin{lemma}
\label{lrot1} Assume that the hypothesis of Lemma
\ref{root3} holds. Let $z$ be the eigenvalue of the
matrix $L$ in the interval $\Delta$ and $f$, $f\neq0$,
the corresponding eigenvector, $Lf=zf$. Then the
(acute) angle $\theta$ between the vectors $f$ and
$f_0=(1,0,0)^\intercal$ satisfies the following bound
\begin{equation}
\label{bo1}
\tan^2\theta\leq\Xi(|\Delta|,d,\|B\|),
\end{equation}
where the function $\Xi$ is given by \eqref{Xi}.
\end{lemma}
\begin{proof}
Assume, without loss of generality, that
$\gamma_+=-\gamma_-=\gamma>0.$
Otherwise one can simply make the corresponding shift of the origin
of the spectral parameter axis. Assume, in addition, that $B\neq0$
and $\lambda\geq0$. (There is no loss of generality in the latter
assumption since, for $\lambda<0$, instead of $L$ one may consider
the matrix $-L$.)

Thus, in the proof we will assume that
$$
\Delta=(-\gamma,\gamma), \quad 0\leq\lambda<\gamma, \text{ \, and \,
} d=\min(\gamma-\lambda,\lambda+\gamma).
$$

Under the hypothesis that
$\|B\|<\sqrt{d|\Delta|}\,\,(=\sqrt{2d\gamma})$, from \cite[Theorem 1
(i)]{KMM3} it follows that if the eigenvalue $z$ of $L$ is in
$\Delta$ then the corresponding eigenvector $f$, $Lf=zf$, may be
chosen in the form
\begin{equation*}
    f=(1,x_-,x_+)^\intercal,
\end{equation*}
with $x_\pm\in\bbC$ such that the matrix $X=(x_- \,\,
x_+)^\intercal$ satisfies the Riccati equation
\begin{equation}\label{RicScal}
\lambda X-A_1X+XBX=B^*.
\end{equation}
Moreover,
\begin{equation}\label{zBX}
z=\lambda+BX.
\end{equation}
Taking into account \eqref{BA1} equations \eqref{RicScal}
and \eqref{zBX} imply
\begin{equation}\label{xx}
    x_-=\dfrac{\,\,\,b_-}{\gamma+z} \text{ \, and \, }
    x_+=\dfrac{\,\,\,b_+}{-\gamma+z}.
\end{equation}
Hence
\begin{equation}\label{nX}
    \|X\|^2=\dfrac{|b_-|^2}{(\gamma+z)^2}+\dfrac{|b_+|^2}{(-\gamma+z)^2}.
\end{equation}
In addition, from \eqref{zBX} and \eqref{xx} one concludes that $z$ is
the solution to equation
\begin{equation}
\label{xxx}
z=\lambda+\dfrac{\,\,\,|b_-|^2}{\gamma+z}+\dfrac{\,\,\,|b_+|^2}{-\gamma+z}.
\end{equation}
Let $t\in[0,1]$ be such that
\begin{equation}
\label{bpa} |b_+|^2=t\|B\|^2
\end{equation}
and, hence,
\begin{equation}
\label{bma}
|b_-|^2=(1-t)\|B\|^2.
\end{equation}
Notice that under the assumptions we use, the bounds
$z_{\mathrm{min}}$ of \eqref{dBm} and $z_{\mathrm{max}}$ of
\eqref{dBp} can  be written in the form
\begin{align}
\label{dBm1}
z_{\mathrm{min}}&=
\dfrac{\gamma+\lambda}{2}-\sqrt{\dfrac{(\gamma-\lambda)^2}{4}+\|B\|^2},\\
\label{dBp1}
z_{\mathrm{max}}&=
-\dfrac{\gamma-\lambda}{2}+\sqrt{\dfrac{(\gamma+\lambda)^2}{4}+\|B\|^2}.
\end{align}
It is easy to see that, given the value of $\|B\|$, for $t$ in
\eqref{bpa} and \eqref{bma} varying between 0 and 1 the solution $z$
to equation \eqref{xxx} fills the whole interval
$[z_{\mathrm{min}},z_{\mathrm{max}}]$. Moreover, with $t$
decreasing from 1 to 0 the value of $z$ is continuously and
monotonously  increasing from $z_{\mathrm{min}}$ to
$z_{\mathrm{max}}$.

On the other hand one can express $t$ through $z$. With $|b_\pm|$
given by \eqref{bpa} and \eqref{bma} from \eqref{xxx} it follows
that
\begin{equation}
\label{alpha}
t=\dfrac{1}{2\gamma\|B\|^2}[(z-\lambda)(z^2-\gamma^2)-\|B\|^2(z-\gamma)].
\end{equation}
Taking this into account, we rewrite expression \eqref{nX}
in the form
$$
\|X\|^2=\varphi(z),
$$
where the function $\varphi$ is given by
\begin{equation}
\label{X2z} \varphi(z)=\dfrac{\|B\|^2+2(\lambda-z)z}{\gamma^2-z^2}.
\end{equation}
That is, given the value of $\|B\|$, the norm of the
solution $X$ to the Riccati equation \eqref{RicScal}
may be considered as a function of the only variable
$z$ that runs through the interval
$[z_{\mathrm{min}},z_{\mathrm{max}}]$.

There is a single point $z_0$ within the interval $(-\gamma,\gamma)$
where the derivative of the function $\varphi(z)$ is zero. This
point reads
\begin{equation}
\label{z0}
z_0=\left\{\begin{array}{cl}
0 & \text{if}\quad \lambda=0,\\[1mm]
\dfrac{2\gamma^2-\|B\|^2}{2\lambda}-
\sqrt{\left(\dfrac{2\gamma^2-\|B\|^2}{2\lambda}\right)^2 -\gamma^2} &
\text{if}\quad \lambda>0.
\end{array}\right.
\end{equation}
It provides the function $\varphi(z)$ with a maximum.

One concludes by inspection that inequality \eqref{Bdd} (along with
the assumptions $\lambda\geq 0$ and $B\neq0$) implies
\begin{equation*}
z_0<z_{\mathrm{max}}.
\end{equation*}
At the same time $z_0 \leq z_{\mathrm{min}}$ if $0<\|B\|\leq \beta$
and $z_0>z_{\mathrm{min}}$ if $\beta<\|B\|<\sqrt{2d\gamma}$ where
\begin{equation}
\label{beta}
\beta=\bigl[(\gamma-\lambda)\sqrt{\gamma}
(\sqrt{\gamma}-\sqrt{\gamma-\lambda})\bigr]^{1/2}=
\left[d\sqrt{|\Delta|}\dfrac{\sqrt{|\Delta|}-\sqrt{2d}}{2}\right]^{1/2}.
\end{equation}
Therefore,
\begin{equation}
\label{fim}
\max\limits_{z\in[z_{\mathrm{min}},z_{\mathrm{max}}]}\varphi(z)
=\varphi(z_{\mathrm{min}})\quad\text{if}\quad 0<\|B\|\leq\beta
\end{equation}
and
$$
\max\limits_{z\in[z_{\mathrm{min}},z_{\mathrm{max}}]}\varphi(z)
=\varphi(z_0)\quad\text{if}\quad \beta<\|B\|<\sqrt{d|\Delta|}.
$$

By substituting \eqref{dBm1} and \eqref{z0} into \eqref{X2z}
one arrives with
\begin{align}
\varphi(z_\mathrm{min})&=\dfrac{d^2}{2\|B\|^2}\left(1+\dfrac{2\|B\|^2}{d^2}-
\sqrt{1+\dfrac{4\|B\|^2}{d^2}}\right)
=\tan^2\biggl(\dfrac{1}{2}\arctan\dfrac{2\|B\|}{d}\biggr)
\end{align}
and
\begin{align}
\varphi(z_0)&=1+\dfrac{2\|B\|^2}{|\Delta|^2}-
\dfrac{2}{|\Delta|^2}\sqrt{
(d|\Delta|-\|B\|^2)\bigl((|\Delta|-d)|\Delta|-\|B\|^2\bigr)},
\end{align}
respectively. To get \eqref{bo1}, it only remains to observe that
$\tan\theta=\|X\|$.

The proof is complete.
\end{proof}

\begin{remark}
\label{rem1} The bound \eqref{bo1} is optimal in the sense that
given the values of $|\Delta|>0$, $d\in(0,|\Delta|/2)$, and
$\|B\|<\sqrt{d|\Delta|}$, it is possible to choose a matrix $L$ of
the form \eqref{Lbl}, \eqref{BA1} such that for the eigenvector
$f=(1,x_-,x_+)^{\intercal}$ associated with the (only) eigenvalue
$z$  of $L$ within the interval $(-\gamma_-,\gamma_+)$ inequality
\eqref{bo1} turns into equality.
\end{remark}

To prove this statement set $\gamma=\frac{|\Delta|}{2}$,
$\gamma_\pm=\pm\gamma$, and $\lambda=\gamma-d$. If $\|B\|\leq\beta$
where $\beta$ is given by \eqref{beta} then choose $b_-=0$ and
$b_+=\|B\|$. Observe that in this case $z=z_{\rm min}$ and hence by
\eqref{fim} such a choice of $b_\pm$ just provides
$\|X\|^2=x_-^2+x_+^2$ with its maximal possible value, i.e. the
equalities $\tan^2\theta=\varphi(z_{\rm min})=\Xi(|\Delta|,d,\|B\|)$
hold. If $\|B\|>\beta$, first compute $t$ by formula \eqref{alpha}
for $z=z_0$ with $z_0$ given by \eqref{z0}. Then introduce
$b_+=\sqrt{t}\|B\|$ and $b_-=\sqrt{1-t}\,\|B\|$. In such a case
$z=z_0$ is the eigenvalue of the matrix $L$ in $\Delta$ and we have
the equality $\tan^2\theta=\varphi(z_0)$, that is, again the equality
$\tan^2\theta=\Xi(|\Delta|,d,\|B\|)$ holds.

\begin{example}
\label{rem2} Again assume that
$\gamma_+=-\gamma_-=\frac{|\Delta|}{2}>0$. Assume in addition that
$\lambda=0$ and $b_+=b_-=\dfrac{b}{\sqrt{2}}$ for some $b\geq 0$.
From \eqref{xxx} it is easy to see that in this case $z=0$ is the
(only) eigenvalue of the matrix $L$ within the interval $\Delta$.
Moreover, for the corresponding eigenvector
$f=(1,x_-,x_+)^\intercal$ by \eqref{xx} one infers that
$x_-=-\dfrac{b}{\sqrt{2}d}$ and $x_+=\dfrac{b}{\sqrt{2}d}$ taking
into account that $\gamma_-=-d$ and $\gamma_+=d$. Since $\|B\|=b$,
the equality $\tan\theta=\sqrt{|x_-|^2+|x_+|^2}$ yields
\begin{equation*}
\tan\theta=\dfrac{\|B\|}{d}.
\end{equation*}
Notice that in this example
$\Xi(|\Delta|,d,\|B\|)=\Xi(2d,d,\|B\|)=\dfrac{\|B\|^2}{d^2}$ and, thus,
the equality $\tan^2\theta=\Xi(|\Delta|,d,\|B\|)$ holds, too.
\end{example}

\begin{example}
\label{remdel} Consider a matrix $L$ of the form \eqref{Lbl} with $\gamma_-$,
$\gamma_+$, and $\lambda$ like in Example \ref{rem2}, that is,
with \mbox{$\gamma_+=-\gamma_-=d>0$} and $\lambda=0$. Set $b_+=0$ and let
$b_-$ satisfy inequalities $0\leq
b_-<\sqrt{d|\Delta|}$. Obviously, $\|V\|=b_-$, $|\Delta|=2d$ and,
thus, we have $\|V\|<\sqrt{2}d$. The eigenvalue
$z$ of the matrix $L$ in the interval $\Delta$ (which is the
corresponding solution to \eqref{xxx}) simply coincides with $z_{\rm
max}$ (cf. formula \eqref{dBp1}),
$$
z=-\dfrac{d}{2}+\sqrt{\dfrac{d^2}{4}+\|V\|^2}
$$
Clearly, $z\to d$ as $\|V\|\to\sqrt{2}d$. That is, in this case the
distance $\delta=\dist(\sigma'_0,\sigma_1)$ between the perturbed
spectral set $\sigma'_0=\{z\}$ and unperturbed spectral set
$\sigma_1=\{-d,d\}$ can be done arbitrarily small.
\end{example}

\section{General case}
\label{Sgener}

Recall that by a finite spectral gap of a self-adjoint operator $T$
one understands an \textit{open} finite interval on $\bbR$ lying
in the resolvent set of $T$ and being such that both of its end
points belong to the spectrum of $T$.

In the sequel, we adopt the following hypothesis.

\begin{hypothesis}
\label{Hyp1}
Let the Hilbert space $\fH$ be decomposed into the orthogonal
sum of two subspaces, i.e.
\begin{equation}
\label{decomp}
\fH=\fH_0\oplus\fH_1.
\end{equation}
Assume that a self-adjoint operator $L$ on $\fH$ reads with respect
to the decomposition \eqref{decomp} as a $2\times2$ operator block
matrix
\begin{equation*}
L=\begin{pmatrix}
  A_0 & B \\
  B^* & A_1
\end{pmatrix}, \qquad \Dom(L)=\fH_0\oplus\Dom(A_1),
\end{equation*}
where $A_0$ is a bounded self-adjoint operator on $\fH_0$, $A_1$ a
possibly unbounded self-adjoint operator on $\fH_1$, and $B$ a
bounded operator from $\fH_1$ to $\fH_0$. Assume in addition, that
$A_1$ has a finite spectral gap $\Delta=(\gamma_-,\gamma_+)$,
$\gamma_-<\gamma_+$, the spectrum of $A_0$ lies in $\Delta$, i.e.
$\spec(A_0)\subset\Delta$, and
\begin{equation}
\label{BdD}
\|B\|<\sqrt{d|\Delta|},
\end{equation}
where
$$
d=\dist(\spec(A_0),\spec(A_1)).
$$
\end{hypothesis}

If $f$ is a non-zero element of the Hilbert space $\fH$ and $\fK$ is
a subspace of  $\fH$, by the angle between $f$ and $\fK$ we
understand the acute angle  $\theta$ between $f$ and its orthogonal
projection $f_\fK$ onto $\fK$, that is,
$\theta=\arccos(\|f_\fK\|/\|f\|)$.

\begin{theorem}
\label{ThMain}
Assume Hypothesis \ref{Hyp1}. Assume in addition that the operator $L$
has an eigenvalue lying in the gap $\Delta$. Let $f$ be an eigenvector
of $L$ associated with this eigenvalue. Then the (acute) angle $\theta$ between
the vector $f$ and the subspace $\fH_0$ satisfies the bound
\begin{equation}\label{bo2}
\tan^2\theta\leq\Xi(|\Delta|,d,\|B\|),
\end{equation}
where the function $\Xi$ is given by \eqref{Xi}.
\end{theorem}
\begin{proof}
Assume that the eigenvector $f=f_0\oplus f_1$, $f_0\in\fH_0$,
$f_1\in\dom(A_1)$, of the operator $L$ is associated with an
eigenvalue $z\in\Delta$. Then the following equalities hold
\begin{align}
\label{e1}
A_0 f_0 + Bf_1 =zf_0&\\
\label{e2} B^* f_0 + A_1 f_1 =zf_1&
\end{align}
Taking into account that $z$ is in the resolvent set of $A_1$, from
\eqref{e2} it follows that
\begin{equation}
\label{e3} f_1=-(A_1-z)^{-1}B^*f_0.
\end{equation}
Hence, $f_0\neq0$ (otherwise, for $f_0=0$, one would have $f_1=0$
and then $f=0$). Equations \eqref{e1} and \eqref{e3} yield
\begin{equation*}
A_0 f_0-B(A_1-z)^{-1}B^*f_0=zf_0,
\end{equation*}
which implies
\begin{equation}
\label{M1} \lal A_0 f_0,f_0\ral-\lal
B(A_1-z)^{-1}B^*f_0,f_0\ral=z\|f_0\|^2
\end{equation}

From now on suppose that
\begin{equation}
\label{f0} \|f_0\|=1
\end{equation}
and set $\lambda=\lal A_0f_0,f_0\ral$. Clearly,
\begin{equation}
\label{lin} \lambda\in[\inf\spec(A_0),\sup\spec(A_0)].
\end{equation}

By the spectral theorem we have
\begin{equation}
\label{BAB} \lal
B(A_1-z)^{-1}B^*f_0,f_0\ral=\int\limits_{\bbR\setminus(\gamma_-,\gamma_+)}
\dfrac{\lal d\sE_{A_1}(\mu)B^*f_0,B^*f_0\ral}{\mu-z},
\end{equation}
where $\sE_{A_1}(\mu)$, $\mu\in\bbR$, denotes the spectral family of
$A_1$. Let
\begin{equation*}
\Delta_-=(-\infty,\gamma_-]\quad\text{and}\quad\Delta_+=[\gamma_+,\infty).
\end{equation*}
By the mean value theorem there are real numbers
$\mu_-\leq\gamma_-$ and $\mu_+\geq\gamma_+$ such that
\begin{align}
\label{E1}
\int_{\Delta_\pm} \dfrac{\lal
d\sE_{A_1}(\mu)B^*f_0,B^*f_0\ral}{\mu-z}&= \dfrac{\lal
\sE_{A_1}\bigl(\Delta_\pm\bigr)B^*f_0,B^*f_0\ral}{\mu_\pm -z}
=\dfrac{\|
\sE_{A_1}\bigl(\Delta_\pm\bigr)B^*f_0\|^2}{\mu_\pm -z},\quad
\end{align}
respectively. Introduce the non-negative numbers $b_\pm$ by
\begin{equation}
\label{bpm}
b_\pm=\sqrt{\alpha_\pm}\|\sE_{A_1}\bigl(\Delta_\pm\bigr)B^*f_0\|,
\end{equation}
where
\begin{equation}
\label{apm} \alpha_\pm=\frac{|\gamma_\pm-z|}{|\mu_\pm-z|}\leq 1.
\end{equation}
Obviously,
\begin{equation}
\label{Ibpm} \int_{\Delta_\pm} \dfrac{\lal
d\sE_{A_1}(\mu)B^*f_0,B^*f_0\ral}{\mu-z}=\dfrac{b_\pm^2}{\gamma_\pm -z}.
\end{equation}
Thus, taking into account \eqref{lin}, \eqref{BAB}, and \eqref{E1},
equation \eqref{M1} turns into
\begin{equation}
\label{Mlam} \lambda - \dfrac{b_-^2}{\gamma_- -z} -
\dfrac{b_+^2}{\gamma_+ -z}=0
\end{equation}
At the same time, by \eqref{e3} we have
\begin{equation}
\label{f1norm}
\|f_1\|^2=\int\limits_{\bbR\setminus(\gamma_-,\gamma_+)} \dfrac{\lal
d\sE_{A_1}(\mu)B^*f_0,B^*f_0\ral}{(\mu-z)^2}.
\end{equation}
The contributions of the intervals $(-\infty,\gamma_-]$ and
$[\gamma_+,\infty)$ to the integral on the r.h.s. part of
\eqref{f1norm} are estimated separately. For the first interval one
derives
\begin{equation*}
\int_{\Delta_-} \dfrac{\lal
d\sE_{A_1}(\mu)B^*f_0,B^*f_0\ral}{(\mu-z)^2}\leq \dfrac{1}{z-\gamma_-}
\int_{\Delta_-} \dfrac{\lal
d\sE_{A_1}(\mu)B^*f_0,B^*f_0\ral}{z-\mu},
\end{equation*}
which by \eqref{Ibpm} means
\begin{equation}
\label{Ibmr} \int_{\Delta_-} \dfrac{\lal
d\sE_{A_1}(\mu)B^*f_0,B^*f_0\ral}{(\mu-z)^2}\leq
\dfrac{b_-^2}{(\gamma_- -z)^2}.
\end{equation}
In a similar way one concludes that
\begin{equation}
\label{Ibpr} \int_{\Delta_+} \dfrac{\lal
d\sE_{A_1}(\mu)B^*f_0,B^*f_0\ral}{(\mu-z)^2}\leq
\dfrac{b_+^2}{(\gamma_+ -z)^2}.
\end{equation}
Then by combining \eqref{f1norm}, \eqref{Ibmr}, and \eqref{Ibpr} one
infers that
\begin{equation}
\label{fxx}
\|f_1\|^2\leq x_-^2+x_+^2,
\end{equation}
where
\begin{equation}
\label{xpm} x_\pm=-\dfrac{b_\pm}{\gamma_\pm-z}.
\end{equation}
From \eqref{Mlam}, \eqref{xpm} it follows that the vector
$y=(1,x_-,x_+)^\intercal$ is an eigenvector of the $3\times 3$
matrix
\begin{equation}
\nonumber
\widetilde{L}=\left(\begin{array}{lll}
\lambda & b_- & b_+ \\
b_- &  \gamma_- & 0 \\
b_+ &  0 & \gamma_+
\end{array}\right)
\end{equation}
associated with the eigenvalue $z$, that is, $\widetilde{L}y=zy$. By
\eqref{lin} for $\delta=\dist(\lambda,\{\gamma_-,\gamma_+\})$ we have
\begin{equation}
\label{deld}
d\leq\delta\leq\dfrac{|\Delta|}{2}.
\end{equation}
In addition, by \eqref{bpm} the square of the norm
$\|\widetilde{B}\|=\sqrt{b_-^2+b_+^2}$ of the $1\times 2$ matrix-row
$\widetilde{B}=(b_-\quad b_+)$ reads
\begin{align}
\nonumber \|\widetilde{B}\|^2&=\alpha_-^2\lal
\sE_{A_1}\bigl(\Delta_-\bigr)B^*f_0,B^*f_0\ral+
\alpha_+^2\lal \sE_{A_1}\bigl(\Delta_+\bigr)B^*f_0,B^*f_0\ral
\end{align}
and hence
\begin{align}
\nonumber
\|\widetilde{B}\|^2\leq &\lal
\sE_{A_1}\bigl(\Delta_-\bigr)B^*f_0,B^*f_0\ral+ \lal
\sE_{A_1}\bigl(\Delta_+\bigr)B^*f_0,B^*f_0\ral\\
\nonumber
&\quad=\lal B^*f_0,B^*f_0\ral=\|B^*f_0\|^2\\
\label{BtB}
\leq & \|B\|^2,
\end{align}
taking into account first \eqref{apm} and then \eqref{f0}. By the
hypothesis inequality \eqref{BdD} holds. Combining \eqref{deld} and
\eqref{BtB} with \eqref{BdD} implies
\begin{equation}
\|\widetilde{B}\|^2<\sqrt{\delta|\Delta|}.
\end{equation}
By Lemma \ref{lrot1} one then concludes that
$x_-^2+x_+^2\leq\Xi(|\Delta|,\delta,\|\widetilde{B}\|)$ which
by \eqref{f0} and \eqref{fxx} implies that
\begin{equation}
\label{teprel}
\tan^2\theta\leq\Xi(|\Delta|,\delta,\|\widetilde{B}\|).
\end{equation}
Given $|\Delta|>0$, $d\in(0,|\Delta|/2]$, and $\|B\|$ satisfying
\eqref{BdD}, it is easy to see that the function
$\Xi(|\Delta|,\delta,\|\widetilde{B}\|)$ is monotonously increasing
with increasing $\|\widetilde{B}\|$, $\|\widetilde{B}\|\leq\|B\|$.
For $d<|\Delta|/2$ it also monotonously increases if $\delta$
decreases from $\frac{|\Delta|}{2}$ to $d$. Therefore, from
\eqref{teprel} it follows that
$\tan^2\theta\leq\Xi(|\Delta|,d,\|B\|)$, completing the proof.
\end{proof}

\begin{remark}
\label{rem1a} The bound \eqref{bo2} is optimal. This follows
from Remark \ref{rem1}.
\end{remark}

\begin{remark}
\label{remms} Notice that under condition
$\|B\|<\sqrt{d(|\Delta|-d)}$ by \cite[Theorem
5.3]{MS01} the
operator angle $\Theta$ between the unperturbed and
perturbed spectral subspaces $\Ran\sE_A(\sigma_0)$ and
$\Ran\sE_L(\sigma'_0)$ satisfies the following (sharp) estimate:
\begin{equation}
\label{thms}
\Theta\leq\frac{1}{2}\arctan\kappa(\|B\|),
\end{equation}
where the function $\kappa(b)$ is defined for $0\leq b<\sqrt{d(|\Delta|-d)}$ by
\begin{equation*}
\kappa(b)=\left\{\begin{array}{ll} \displaystyle \frac{2b}{d}
& \text{\, if \,} b \le
\displaystyle\sqrt{\frac{d}{2}\left(\frac{|\Delta|}{2}-d\right)},\\
\displaystyle \frac{b\dfrac{|\Delta|}{2} + \sqrt{d(|\Delta| -
d )\Bigl[ \Bigl(\dfrac{|\Delta|}{2} -  d\Bigr)^2 + b^2
\Bigr]}} {d(|\Delta|-d) - b^2} & \text{\, if \,} b
>\displaystyle\sqrt{\frac{d}{2}\left(\frac{|\Delta|}{2}-d\right)}.
\end{array}\right.
\end{equation*}
Surely, the bound \eqref{thms} implies the corresponding estimate
for the angle $\theta$:
\begin{equation}
\label{thms1}
\theta\leq\dfrac{1}{2}\arctan\kappa(\|B\|)\quad\text{whenever}\quad
\|B\|<\sqrt{d(|\Delta|-d)}.
\end{equation}
One observes by inspection that $\Xi(|\Delta|,d,b)\leq
\tan^2\left(\frac{1}{2}\arctan\kappa(b)\right)$, $0\leq
b<\sqrt{d(|\Delta|-d)}$. Moreover, if $|\Delta|>2d$
then for
$\sqrt{\frac{d}{2}\left(\frac{|\Delta|}{2}-d\right)}<b<\sqrt{d(|\Delta|-d)}$
the strict inequality
$\Xi(|\Delta|,d,b)<\tan^2\left(\frac{1}{2}\arctan\kappa(b)\right)$
holds. Therefore, the bound \eqref{thms1} is not optimal in the case
of eigenvectors.
\end{remark}

Now we are in position to prove Theorem \ref{Th1}. This theorem
appears to be a simple corollary to Theorem \ref{ThMain}.

\begin{proof}[Proof of Theorem \ref{Th1}]
Set $\fH_0=\Ran\sE_A(\sigma_0)$ and $\fH_0=\Ran\sE_A(\sigma_1)$.
With respect to the orthogonal
decomposition $\fH=\fH_0\oplus\fH_1$ the operators $A$ and $V$ read
as $2\times2$ block operator matrices,
\begin{equation*}
A=\begin{pmatrix}
  A_0 & 0 \\
  0 & A_1
\end{pmatrix} \quad
\text{and} \quad
V=\begin{pmatrix}
  0 & B \\
  B^* & 0
\end{pmatrix},
\end{equation*}
where $B=V|_{\fH_1}$; \,\, $\Dom(A)=\fH_0\oplus\Dom(A_1)$ and
$\Dom(L)=\Dom(A)$. Assume that $\Delta$ is a gap of the set
$\sigma_1$ that contains the whole set $\sigma_0$. Surely, the
length $|\Delta|$ of this gap satisfies the estimate
$|\Delta|\geq 2d$ and the bound \eqref{Vd2} implies
the inequality $\|B\|<\sqrt{d|\Delta|}$. Then by Theorem \ref{ThMain}
we have
\begin{equation}
\nonumber
\tan^2\theta\leq\Xi(|\Delta|,d,\|V\|),
\end{equation}
taking into account that $\|V\|=\|B\|$. Now it only remains to
observe that $\Xi(D,d,\|V\|)$ is a non-increasing function of the
variable $D$, $D\geq 2d$. For $D$ varying in the interval
$[2d,\infty)$ it achieves its maximal value just at $D=2d$ and this
value equals
\begin{equation}
\nonumber
\max\limits_{D: \,D\geq 2d}\Xi(|\Delta|,d,\|V\|)=\dfrac{\|V\|^2}{d^2}.
\end{equation}
Thus, the following inequality holds
$$
\tan\theta\leq\dfrac{\|V\|}{d}.
$$
The proof is complete.
\end{proof}
\begin{remark}
Example \ref{rem2} shows that the bound \eqref{thet} is sharp.
\end{remark}
\subsection*{Acknowledgment}
We greatfully acknowledge that this work was supported
by the Deutsche Forschungsgemeinschaft (DFG), the
Heisenberg-Landau Program, and the Russian Foundation
for Basic Research.

\end{document}